# Stochastic Multi-class Traffic Assignment for Autonomous and Regular Vehicles in a Transportation Network

**S. Roozbeh Mousavi**
Department of Civil Engineering
Sharif University of Technology, Tehran, Iran, 1458889694
Email: roozbeh.mousavi@sharif.edu
0000-0003-4429-7224

**Alireza Yazdiani***
Department of Civil Engineering
Sharif University of Technology, Tehran, Iran, 1458889694
Email: alireza.yazdiani@sharif.edu
0000-0002-4135-5209

**Yousef Shafahi**
Department of Civil Engineering
Sharif University of Technology, Tehran, Iran, 1458889694
Email: shafahi@sharif.edu
0000-0003-4267-4348

* Corresponding author

**ABSTRACT**

A transition period from regular vehicles (RVs) to autonomous vehicles (AVs) is imperative. This article explores both types of vehicles using a route choice model, formulated as a stochastic multi-class traffic assignment (SMTA) problem. In RVs, cross-nested logit (CNL) models are used since drivers do not have complete information and the unique characteristics of CNL. AVs, however, are considered to behave in a user equilibrium (UE) due to complete information about the network. The main innovation of this article includes the introduction of three solution methods for SMTA. Depending on the size of the network, each method can be used. These methods include solving the nonlinear complementary problem (NCP) with GAMS software, the decomposition-assignment algorithm, and the modified Wang's algorithm. Through the modification of Wang's algorithm, we have increased the convergence speed of Wang's algorithm and shown its numerical results for the Nguyen and Sioux Falls networks. As it is not possible to consider all paths in the traffic assignment, we proposed a creative path generation-assignment (PGA) algorithm. This algorithm generates several attractive paths for each origin-destination (OD), and the modified Wang's algorithm assigns traffic flow.

## INTRODUCTION

As AVs enter the transportation network, traffic parameters such as travel time, mode choice models, and route choice models will also change. Compared to RVs, these vehicles will be able to drive at a closer distance to the one in front. We will also see an increase in the capacity of the transportation network and its safety due to the exchange of information between the vehicles and the infrastructure network. This transition period from RVs to AVs will be a relatively significant one in terms of duration, therefore it is necessary to consider both phenomena simultaneously. In this research, the special characteristics of each vehicle are presented, such as (headway, connections with other cars, fuel consumption, link capacity for each vehicle, etc.), as well as the attitudes of the drivers (such as VOT, the comfort of travel, etc.).

Traffic assignment based on the UE approach is one of the oldest approaches to solving the traffic assignment problem. The main assumption of this approach is that all drivers have complete information about the network. In RVs, however, the behavior of drivers is based on a stochastic approach. As a result, the stochastic user equilibrium (SUE) approach has been used to model RVs in many studies. As for AVs, however, since they use advanced technology, it is more appropriate to assume their behavior based on the UE approach. In this research, we will focus on the modeling of SMTA for AVs and RVs, where AVs will behave according to the UE approach and RVs will follow the SUE approach, and specifically for this type, we will use the CNL model.

With the entry of AVs into the network, the network's capacity will increase; in the multi-class mode, there will be more traffic in the links than for RVs. We also expect that based on the characteristics of the CNL model, the flow of RVs in the network will be closer to reality since this model can distinguish between paths that overlap in one or more links and routes that do not overlap. The study assumes that all links can be used by both RVs and AVs. As previously mentioned, RV users will behave according to the CNL model when choosing a route, whereas AV users will follow the UE model. In order to determine the travel cost for links, it is only necessary to know the travel time and the fuel cost of the links.

This study aims to introduce a suitable algorithm for solving the SMTA model for AVs and RVs, due to the nonlinear and complex structure of the CNL model. As a general rule, three solution methods have been introduced in this research. All of the algorithms are path-based, and each was designed with a specific purpose in mind. Their accuracy and convergence speed will be evaluated.

The first solution method involves implementing NCP in the GAMS software environment and solving the problem with solver Konopt4. The accuracy of the final answer is one of the benefits of this method. This method also has disadvantages, such as computational limitations in GAMS software, which make it unsuitable for medium-sized and large networks. The second method uses the decomposition algorithm to decompose the transportation network into OD units, which can be solved using GAMS software. For each iteration of this method, traffic assignments will be performed locally for each OD. GAMS software is used to solve the traffic assignment problem for each OD in this method. As a final step, we will present a modified version of Wang's algorithm. Through slight modifications to Wang's algorithm, this algorithm will significantly increase the speed of its convergence, as shown in the numerical example.

The three methods proposed in this study, and most path-based traffic assignment algorithms, require initial paths in order to function. Therefore, selecting the right path is crucial in determining the accuracy of the answer and its deviation from the state when all paths are considered. This research presents a PGA algorithm as a solution to this problem. The traffic assignment and path generation processes will be performed alternately in this algorithm. In this algorithm, Yen's k-Shortest Path algorithm is used to generate paths, and a modified Wang's algorithm is used to assign traffic flow.

## LITERATURE REVIEW

Considering that this article focuses on AVs, discrete choice models, and traffic assignment in the transportation network, we will have a brief overview of previous research in these areas as follows:

### Autonomous vehicles

In comparison to RVs, AVs have unique and distinctive capabilities that will have an impressive impact on traffic parameters and driving behavior. In addition to increased network capacity and safety, driverless operation, fuel consumption reduction, and environmental benefits, these vehicles are cost-effective and environmentally friendly. Vehicle fuel consumption will be partially reduced due to the optimal behavior of AVs and their considerable impact on reducing traffic congestion. Unlike RVs, these vehicles do not require a driver, so their users can engage in other activities while traveling, which decreases the value of time (VOT) for drivers [1].

### Discrete choice models

The logit and probit models are well-known discrete choice models. Polynomial logit models have closed-form probability functions and are computationally efficient when used for traffic assignments. For modeling route choice behavior, the logit model suffers from the overlapping problem and considers the same probability distribution for all error terms in the utility function. Alternatively, the polynomial probit model considers the correlation between paths and topographic features in the network. Despite the theoretical advantages of this model, it lacks a closed-form probability function and is computationally very expensive for traffic assignment.

In order to improve and eliminate the shortcomings of the polynomial logit model, many models have been proposed. Among these models are the C-logit model [2]and the path-size logit model [3], [4]. In these two models, the systematic or observed term of the utility function is modified to consider the similarities and correlation between the possible paths of each OD. The generalized extreme value theory (GEV) is another class of these models. Initially, McFadden [5] proposed the GEV models class. Among these are the paired combinatorial logit (PCL) model [6], [7], the kernel logit model [8], and the cross-nested logit (CNL) model [9]. Models such as these attempt to model the correlation between paths in more realistic ways by considering more comprehensive structures. Vovsha [10] used the CNL model for mode choice in the third stage of the four-level model. Eventually, Prashker and Bekhor [6] and Vovsha and Bekhor [9] utilized this model for route choice.

### Multi-class and stochastic traffic assignment

In the definitive approach, users' route choices are solely influenced by travel costs. As a result, the UE problem will become the Wardrop user equilibrium problem. Wardrop [11] first presented this law in the following form:

> **Wardrop User Equilibrium law:** In equilibrium conditions, users cannot improve their travel costs by unilaterally changing their routes.

Wardrop's user equilibrium law can also be formulated stochastically, in the same way as the definitive user equilibrium:

> **Stochastic user equilibrium (SUE)**: In equilibrium conditions, users cannot improve their <u>perceived</u> travel cost by unilaterally changing their route.

One disadvantage of primary UE models is that they assume all users have the same preferences and decision metrics. The UE model is rewritten by Dafermos [12] as a group of users, assuming that the link travel cost varies from the users' perspective and that each user has their preferences. As a result of this model, every OD's demand is divided into different groups. In the route choice stage, each of these groups has the same cost function for the link, so the cost of travel on each path will be the same. The model introduced by Dafermos is known as the multi-class user equilibrium (MUE) model.

Ying and Wang, 2018 [13] considered multi-class flows and developed an optimal network design by establishing special routes for AVs and RVs and adjusting pricing for the AVs. In 2016, Chen et al. [14]examined how to allocate lines to AVs based on homogeneous multi-class flows and AV penetration rate into the transportation network. In 2017, Ghiasi et al. [15] examined and analyzed issues such as the capacity and dedication of specific routes to AVs, considering the multi-class flow.

In 2019, Wang et al. [16] modeled the network traffic assignment problem in the context of multi-class flows of AVs and RVs as a variational inequality. In addition, they provide recommendations for improving network performance by analyzing effective parameters like demand and link capacity.

To solve the problem of SMTA, several algorithms have been proposed. Most of these algorithms are based on variational inequality (VI) models. Among these algorithms are the projection method [17]–[19], the Tikhonov regularization method [20], and the proximal point method [21]. Each iteration of the mentioned algorithms solves a separate sub-problem to determine the descent direction. Due to this problem, the above algorithms are not suitable for solving the problem of SMTA. Traffic assignment models presented in this study have complex travel cost functions, and solving repeated problems will be very expensive because of the complex expressions.

Flow swapping algorithms between paths are algorithms that show the process of formation of balanced flows in the network based on the behavior of network users. A balanced solution is reached by moving from a non-equilibrium solution and exchanging the flow between paths. Algorithms such as these mimic network users who move toward forming an equilibrium flow by simulating their behavior.

Despite the relative improvement of the flow swapping algorithm, many studies indicate the low convergence speed of their algorithm [16]. Wang et al. Presented a creative algorithm to solve the problems. They use a function called $\emptyset$ to calculate the descent direction. This function swaps the flow between all possible paths within each OD. Swapping flows between routes are also affected by the flow in each route and the cost difference between each route. As part of their algorithm, they have also developed an innovative method for determining the step size to guarantee the feasibility of the solution. In comparison with previous similar algorithms, Wang et al. Claim their algorithm will have a high convergence speed. Their innovative algorithm is called RSRS-MSRA. We will refer to their algorithm as Wang's algorithm in this paper.

## MATHEMATICAL MODEL

RVs and AVs have fundamentally different structures and behaviors. AVs use vehicle-to-vehicle (V2V) communication equipment and vehicle-to-infrastructure (V2I) communication equipment. However, in RVs, all navigation operations fall under the driver's responsibility. As a means of considering these two phenomena together, we will use the SMTA as our methodology. Wang et al. [16] state three main reasons for using the SMTA:

1- RV users' route choices can be significantly influenced by the flow of AVs. As well, link and path travel costs will differ if both types of vehicles are present versus only one type. In addition, the cost difference between the two types of vehicles is due to behavioral asymmetry.
2- From the perspective of the network users, path costs are determined by several factors. Fuel costs and path travel time are the most important elements. Path travel cost can be calculated using VOT, which is a parameter that converts path travel time to path travel cost; The VOT is more valuable to RV drivers than AV drivers since AVs cannot be controlled by the user, so they can engage in other activities.
Various policies have been proposed to promote the use of AVs. For example, AVs can be given their own lanes and paths. In this case, the AVs will be able to use a special path that is specially designed for them. Therefore, SMTA models must be used to consider this distinction.

Consider the transportation network $G\ (\boldsymbol{N}, \boldsymbol{\Gamma})$. $\boldsymbol{N}$ represents the set of nodes, and $\boldsymbol{\Gamma}$ represents the set of links in the network. This network includes RVs and AVs. Therefore, in the formulas of this study, $R$ represents RVs, and $A$ represents AVs.

**Link travel cost function**

The Link travel cost consists of two parts: $\bar{t}_a$, the travel time, and $E_a$, the fuel consumption cost. As mentioned before, due to the difference in the navigation structure of the two types of vehicles, the link's travel time depends on the percentage of each vehicle. In order to determine link travel time, we use the BPR function as follows:

$$\bar{t}_a(x_{a,R}, x_{a,A}) = \bar{t}_0\left(1 + 0.15\left(\frac{x_{a,R} + x_{a,A}}{Q_a}\right)^4\right) \qquad \forall a \epsilon \Gamma \qquad (1)$$

In the above relation, $\bar{t}_a$ represents the travel time and $\bar{t}_0$ represents the free travel time of link $a$, respectively. Also, $x_{a,m}$ represents the flow of a vehicle of type $m$ in link $a$, and $Q_a$ represents the mix capacity of link $a$. Levin and Boyles [22], assuming that the reaction times of both types of vehicles are equal, propose the $Q_a$ relation as follows:

$$Q_a = \frac{1}{\frac{x_{a,R}}{x_{a,R} + x_{a,A}}\frac{1}{Q_{a,R}} + \frac{x_{a,A}}{x_{a,R} + x_{a,A}}\frac{1}{Q_{a,A}}} \qquad \forall a \epsilon \Gamma \qquad (2)$$

In the above equation, $Q_{a,R}$ indicates the capacity of link $a$ , in the case that all vehicles are AVs. Similarly, $Q_{a,A}$ indicates the capacity of link $a$ in the case that all vehicles are AVs. Based on previous studies and the use of real data, Levin and Boyles [22] have proposed the following experimental relationship to calculate the amount of fuel consumed in link $a$:

$$E_a = \frac{l_a}{36.44}\left(14.58\left(\frac{l_a}{\bar{t}_a}\right)^{-0.625}\right) \qquad \forall a \epsilon \Gamma \qquad (3)$$

In the above equation, $l_a$ is equal to the length of link a. Link travel cost is calculated as follows:

$$t_{a,m} = \bar{t}_a(x_{a,R}, x_{a,A}).VOT_m + \eta.E_a \qquad \forall a \epsilon \Gamma, m \epsilon \{RV, AV\} \qquad (4)$$

The $VOT_m$ in the above relation represents the VOT for the user of vehicle type $m$. Also, $\eta$ represents the price of a gallon of gasoline, which in this study is estimated at \$5.5.

**Path travel cost**

Travel costs for each path can be calculated by adding the travel costs of the links associated with that path.

$$c_{k,m}^{w} = \sum_a t_{a,m} . \delta_{a,k,m}^{w} \qquad \forall w, k, m \qquad (5)$$

For the vehicle type $m$, $\delta^{w}_{a,k,m}$ is one if the link a is in the path $k$ between the OD pair $w$; otherwise, it is 0. The path perceived cost for RV users was shown by Wang et al. [16]as follows:

$$C^{w}_{k,R} = c^{w}_{k,R} - \frac{u}{\theta} H^{w}_{k,R} + \frac{u}{\theta} \ln\left(\frac{f^{w}_{k,R}}{q^{w}_{R}}\right) \tag{6}$$

Where $H^{w}_{k,R}$ is as follows:

$$H^{w}_{k,R} = \ln\left[\sum_{b\in\Gamma^{R}} \left(\alpha^{w}_{b,k}\right)^{1/u} \left(\sum_{l\in\Lambda^{w}_{R}} \left[\alpha^{w}_{b,l}\, exp\left(-\theta c^{w}_{l,R}\right)\right]^{1/\mathrm{u}}\right)^{u-1}\right] \tag{7}$$

$$\alpha_{b,k} = \frac{l_b}{l_k} \delta_{b,k} \tag{8}$$

In **Equation 6,** the variable $f^{w}_{k,R}$ represents the flow of RVs in the path $k$ between the OD pair $w$, $q^{w}_{R}$ represents the demand of the RVs between the OD pair $w$ and $\alpha^{w}_{b,k}$ shows the degree of overlapping, which is calculated by **Equation 8**. Also $\Lambda^{w}_{R}$ represents the set of RV paths between OD pair $w$. In this **Equation 8** $l_b$ is the length of the link $b$ and $l_k$ represents the length of the path $k$. Also, if the path $k$ contains the link $b$, the value of $\delta_{b,k}$ is equal to one; otherwise, it is equal to zero.

Wang et al. [16] have shown that under equilibrium conditions, the flow in the path of RVs will be in accordance with the CNL model if and only if the perceptual travel costs of the paths are equal. The perceptual travel cost of AVs users due to complete information from the network will be equal to the observed travel cost of the path, i.e.

$$C^{w}_{k,A} = c^{w}_{k,A} \tag{9}$$

The model chosen for the multi-class equilibrium flow in this study is a complementarity model as follows:

$$f^{w}_{k,m}\left(C^{w}_{k,m} - u^{w}_{m}\right) = 0 \qquad \forall\, w, k, m \tag{10}$$

$$C^{w}_{k,m} - u^{w}_{m} \geq 0 \qquad \forall\, w, k, m \tag{11}$$

$$\sum_{k} f^{w}_{k,m} = q^{w}_{m} \qquad \forall\, w, m \tag{12}$$

$$f^{w}_{k,m} \geq 0 \qquad \forall\, w, k, m \tag{13}$$

In the above equations, $u^{w}_{m}$ represents the lowest perceived cost between OD pair $w$ by the user of vehicle type $m$. It is also implicitly stated that the flow in the links is as follows:

$$x_{a,m} = \sum_{w}\sum_{k} f^{w}_{k,m} \,.\, \delta^{w}_{a,k,m} \qquad \forall\, a, m \tag{14}$$

## SOLUTION METHODS

This article presents three methods for solving the SMTA problem. These methods have been described below, along with their advantages and disadvantages:

**Solving the complementary model with GAMS**

In this method, the problem of SMTA is implemented and solved as a complementary model in the GAMS environment. This kind of problem whose constraints are nonlinear is called a nonlinear complementary program (NCP). The structure of this method is such that the entire network is considered at once, and the traffic assignment operation is performed simultaneously for all classes of vehicles and all ODs. The noteworthy point is that this method can only be used for small networks due to the computational limitations of the GAMS software, such as the number of allowed variables, the number of allowed constraints, etc.

**Algorithm of decomposition and assignment**

This method is an iterative algorithm in which the transportation network is decomposed into OD units. In each iteration, an NCP is solved for each OD, equilibrium conditions are established locally for each OD, and the network is updated.

Unlike the previous method, this iterative algorithm solves one NCP for each OD in each iteration. It is clear that the number of variables considered for each NCP of this method is much less than the number of variables considered in the previous method. Therefore, this method gives us the possibility to solve the problem of traffic assignment by considering the complementary model for medium and large networks by considering the limitations of the GAMS software, i.e., the number of allowed variables and the number of allowed limitations.

**Modified Wang's algorithm**

The Modified Wang's algorithm used in this research has two advantages over the Wang's algorithm in terms of the direction of decreasing the objective function and movement step. In this study, with changes in the direction of descent and the movement step, the performance of the Wang's algorithm has been significantly improved. The function for the proposed descent direction of the present study is as follows:

$$\phi_{k,m}^{w}(\mathbf{f}(\boldsymbol{n})) = \sum_{g\epsilon\Lambda_m^w} \Big[ f_{g,m}^{w}(n)\left(C_{g,m}^{w}(\mathbf{f}(\boldsymbol{n})) - C_{k,m}^{w}(\mathbf{f}(\boldsymbol{n}))\right)_+^{\mu_m} - f_{k,m}^{w}(n)\left(C_{k,m}^{w}(\mathbf{f}(\boldsymbol{n})) - C_{g,m}^{w}(\mathbf{f}(\boldsymbol{n}))\right)_+^{\mu_m} \Big] \quad \forall w,k,m \tag{15}$$

In the above equation, $\mu_m$ is the flow swapping degree of the vehicle type $m$. The parameter $\mu_m$ indicates the intensity of the swapped flow between the paths according to the type of vehicle $m$. The rest of the variables are already defined. In this algorithm, the movement step of the Wang's algorithm is also changed and rewritten as follows:

$$\beta(n) = \begin{cases} \dfrac{1}{h(n)} * \dfrac{1}{\gamma(n)}, & if\ O(n) > O(n-1) \\ \dfrac{1}{h(n)} * \dfrac{O(n)}{O(n-1)}, & if\ O(n) \leq O(n-1) \end{cases} \tag{16}$$

As can be seen, the movement step in the proposed algorithm of this study is considered as two criteria. **Equation 17** represents the sum of the positive values of the vector components of $\phi$ in the Nth iteration, which measures the distance or proximity to the equilibrium solution. The proposed relation $O(n)$ is as follows:

$$O(n) = \sum_{w,k,m} |\phi_{k,m}^{w}(\mathbf{f}(\boldsymbol{n}))| \tag{17}$$

The function $h(n)$ in **Equation 16** guarantees the possibility of an answer in each iteration. The proposed relation of this research for $h\ (n)$ is as follows:

$$h(n) = max\left\{h_{k,m}^{w}(n) \middle| h_{k,m}^{w}(n) = \frac{-\phi_{k,m}^{w}\big(f(n)\big)}{f_{k,m}^{w}(n)}, \quad \forall w,k,m \middle| \left(f_{k,m}^{w} \neq 0\,, \phi_{k,m}^{w} \leq 0\right)\right\} \tag{18}$$

Details of the steps of the Modified Wang's algorithm are given below:

**Steps of the Modified Wang's algorithm:**

**Step 0** Initial Values: Select a value $G > 0$ for the convergence criterion. Put $n = 1, G(n) = \infty, \gamma(n) = 2, \lambda_1 = 0.0001, \lambda_2 = 0.0001, \mu_A = 1, \mu_H = 0.85$. For each $m$ and $w$, select the set of paths $\Lambda_m^w$. Assign demand values to paths as follows:

$$f_{k,m}^{w}(n) = \frac{q_m^w}{|\Lambda_m^w|} \qquad \forall w,k,m \tag{19}$$

**Step 1** Update the flow in the links according to the below relationships and calculate the general cost values of the link and the path.

$$x_a^m(n) = \sum_{w \epsilon W} \sum_{k \epsilon \Lambda_m^w} f_{k,m}^{w}(n).\, \delta_{a,k,m}^{w} \qquad \forall a,m \tag{20}$$

$$\bar{t}_a^m(n) = t_a\big(x_a^R(n), x_a^A(n)\big)\,.VOT^m + \eta\,.E_a(n) \qquad \forall a,m \tag{21}$$

$$c_{k,m}^{w}(n) = \sum_{a} \bar{t}_a^m(n).\, \delta_{a,k,m}^{w} \qquad \forall\, w,k,m \tag{22}$$

$$C_{k,R}^{w}(n) = c_{k,R}^{w}(n) - \frac{u}{\theta} H_{k,R}^{w}(n) + \frac{u}{\theta} \ln\left(\frac{f_{k,R}^{w}(n)}{q_R^w}\right) \qquad \forall\, w,k \tag{23}$$

$$C_{k,A}^{w}(n) = c_{k,A}^{w}(n) \qquad \forall\, w,k \tag{24}$$

**Step 2** Calculate the following functions:

$$\begin{aligned}\phi_{k,m}^{w}(f(n)) = \sum_{g \epsilon \Lambda_m^w} \Big[ & f_{g,m}^{w}(n)\left(C_{g,m}^{w}(f(n)) - C_{k,m}^{w}(f(n))\right)_+^{\mu_m} \\ & - f_{k,m}^{w}(n)\left(C_{k,m}^{w}(f(n)) - C_{g,m}^{w}(f(n))\right)_+^{\mu_m} \Big] \qquad \forall w,k,m\end{aligned} \tag{25}$$

$$\begin{aligned}h(n) = max\Bigg\{ & h_{k,m}^{w}(n)|\, h_{k,m}^{w}(n) = \frac{-\phi_{k,m}^{w}\big(f(n)\big)}{f_{k,m}^{w}(n)}, \forall w,k,m|(f_{k,m}^{w}(n) \neq 0, \\ & \phi_{k,m}^{w}(\mathrm{f}(n)) \leq 0) \Bigg\}\end{aligned} \tag{26}$$

$$O(n) = \sum_{w,k,m} |\phi_{k,m}^{w}(\mathrm{f}(n))| \tag{27}$$

**Step 3** if $n = 1$, insert $\beta(n) = \frac{1}{h(n)} * \frac{1}{\gamma(n)}$, $\gamma(n) = 9.5$ and go to step 5.

**Step 4** if $O(n) > O(n-1)$, insert $\gamma(n) = \gamma(n-1) + \lambda_1$, $\beta(n) = \frac{1}{h(n)} * \frac{1}{\gamma(n)}$ and go to step 5; otherwise, insert $\gamma(n) = \gamma(n-1) + \lambda_1$, $\beta(n) = \frac{1}{h(n)} * \frac{O(n)}{O(n-1)}$

**Step 5** Update the flow in the paths according to the following equation.

$$f_{k,m}^{w}(n+1) = f_{k,m}^{w}(n) + \beta(n) * \phi_{k,m}^{w}(\mathbf{f}(\boldsymbol{n})) \qquad \forall w,k,m \tag{28}$$

**Step 6** Calculate the value of error and the total cost of travel on the network according to the following equations.

$$G(n) = \frac{\sum_{w,k,m} f_{k,m}^{w}(n)\left(C_{k,m}^{w}(n) - C_{min,m}^{w}(n)\right)}{\sum_{w,k,m} f_{k,m}^{w}(n) C_{k,m}^{w}(n)} \tag{29}$$

$$TC(n) = \sum_{w,k,m} f_{k,m}^{w}(\mathrm{n}) * C_{k,m}^{w} \tag{30}$$

Which $C_{min,m}^{w}$ is calculated from the following equation.

$$C_{min,m}^{w} = \min_{k \epsilon R_m^w}\{ C_{k,m}^{w}(n) \} \tag{31}$$

**Step 7** If $G(n) > G$ put $n = n + 1$ and go to step 1, otherwise stop. An equilibrium answer is obtained.

Comparing the steps of the Modified Wang's and Wang's algorithms, it can be seen that these two algorithms are similar in a large number of steps and differ only in determining the direction of descent and the movement step. It will also be shown that this slight difference will significantly improve the convergence rate.

**Path generation-assignment algorithm**

The previous methods require getting all possible paths for each OD. For small networks, determining all possible paths is not difficult. However, for large and real networks, as well as medium networks such as the Sioux Falls network, which have a lot of ODs and lots of links, it is practically impossible to count all the paths.

The purpose of this paper is to produce a limited number of paths to satisfy two contradictory requests simultaneously. Firstly, the algorithm must have sufficient paths for it to run smoothly. Moreover, an equilibrium solution based on the same limited number of paths is similar to an equilibrium solution based on all paths.

Accordingly, we must answer two fundamental questions in this research based on the above preferences:

1. How did the answer change for the number of different paths for each OD?
2. What is the best path among all possible paths for each OD?

To answer the first question, one needs to perform a relative sensitivity analysis between the number of paths and the algorithm's execution time. On the other hand, it is necessary to create attractive routes from the user's perspective while assigning demand to the network to answer the second question. The algorithm should identify the optimal paths based on the cost of the links and their crowding, as well as the preferences of the users.

In the current research, one of the main assumptions is that RV users do not have complete knowledge of the network situation, such as the travel cost of links, paths, etc., which means that SUE models should be used to model RV users' behavior. To generate the path, it is assumed that there is complete information about the network, and the perceived travel costs are the same as the actual travel costs. Then the route is generated based on that information.

A PGA algorithm is an iterative algorithm that generates paths and assigns traffic alternately, resulting in an equilibrium solution. This algorithm generates attractive paths from the users' perspective in every iteration, and the Modified Wang's algorithm assigns demand values. The operation will continue until the set of produced paths in the new iteration does not differ from the set of paths in the previous iteration or differs by just a tiny amount.

For each OD, Yen's algorithm is used to generate attractive paths from the perspective of network users. Following is a description of the steps in the PGA algorithm.

**Steps of path generation- assignment algorithm**

**Step 0** Initial values: Choose a $G > 0$ value for the convergence criterion, an $E > 0$ value for the replication error criterion, and a $k$ value for the number of paths for each OD.
Put: $m = 1, E(m) = \infty, G(m) = 0.1$ .

**Step 1** For each link and both types of vehicles, calculate its free travel cost.

**Step 2** Using Yen's algorithm and according to the links travel cost, for each OD and both types of vehicles, generate the k-shortest path, then update the set of network paths.

**Step 3** Using the Modified Wang's algorithm and considering the convergence criterion $G(m)$, assign the demand values to the network and then update the travel cost of the links and $TC(m)$, the total travel cost of the network.

**Step 4** If $m < 2$, set $m = m + 1$, then go to step 2. Otherwise, calculate the iteration error value according to the following equation.

$$E(m) = \frac{TC(m) - TC(m-1)}{TC(m)} \tag{32}$$

**Step 5** If $E(m) > E$, put $m = m + 1$ and go to step 2; otherwise, stop and put: $m = m + 1$. Then calculate the equilibrium answer using the Modified Wang's algorithm and consider the $G$ convergence criterion.

## NUMERICAL EXAMPLES

In this section, we present the results of implementing the SMTA solution methods for two networks, Nguyen and Sioux Falls.

### Nguyen network

Nguyen's network [23] is given in **Figure 1**. For predetermined demand, four traffic assignment solution methods have been implemented. In all these methods, paths are provided. **Table 1** summarizes the results of implementing all four methods. **Table 1** shows that method one has higher accuracy than all other methods and has the shortest execution time. Based on the results, Compared with the Wang's algorithm, the Modified Wang's algorithm performs eight times better in this network.

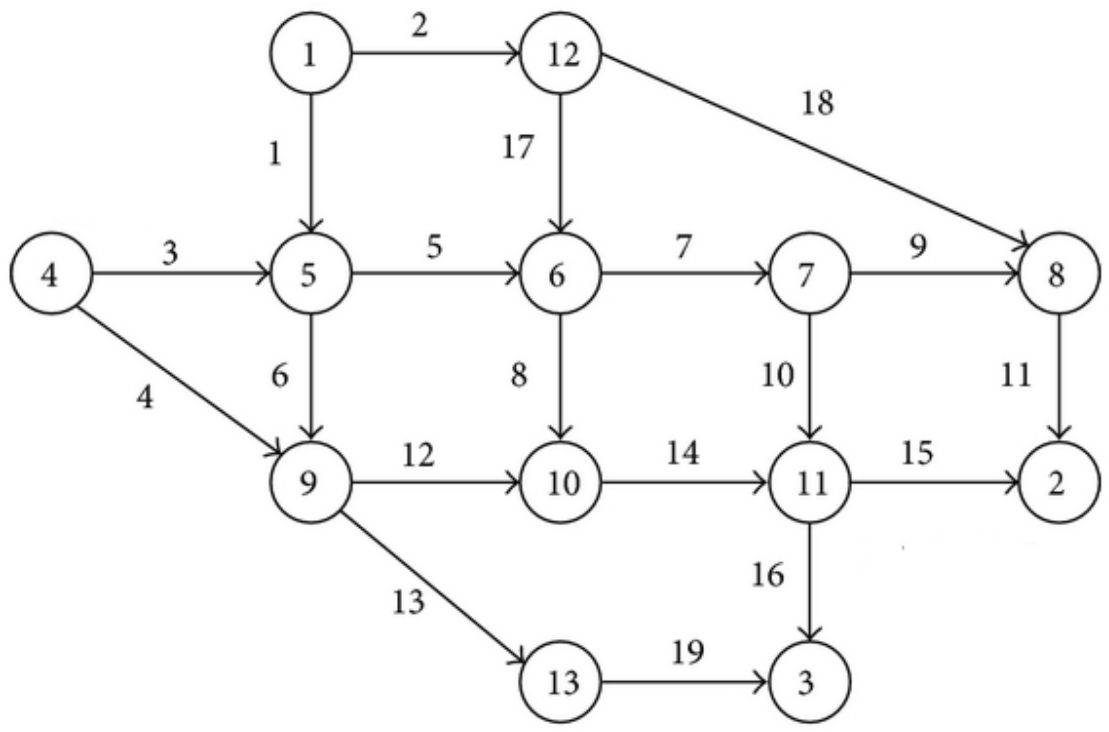

**Figure 1 Nguyen network.**

**Table1 Four traffic assignment methods for the Nguyen network.**

| Method number | solution method | Error | number of iterations | Execution time (seconds) | Total travel cost ($) |
|---|---|---|---|---|---|
| 1 | NCP with GAMS | 6.5E-09 | ----- | 0.015 | 31420.52 |
| 2 | Decomposition algorithm | 0.0001 | 16 | 0.890 | 31421.43 |
| 3 | Wang's algorithm | 0.0001 | 1702 | 2.953 | 31419.51 |
| 4 | Modified Wang's algorithm | 0.0001 | 179 | 0.36 | 31419.69 |

In **Figure 2**, the Wang's and Modified Wang's algorithms are compared regarding their convergence process for 1800 iterations. According to this figure, the Modified Wang's algorithm achieves an error value $10^8$ times lower higher than the Wang's algorithm for 1000 iterations, which demonstrates its superiority.

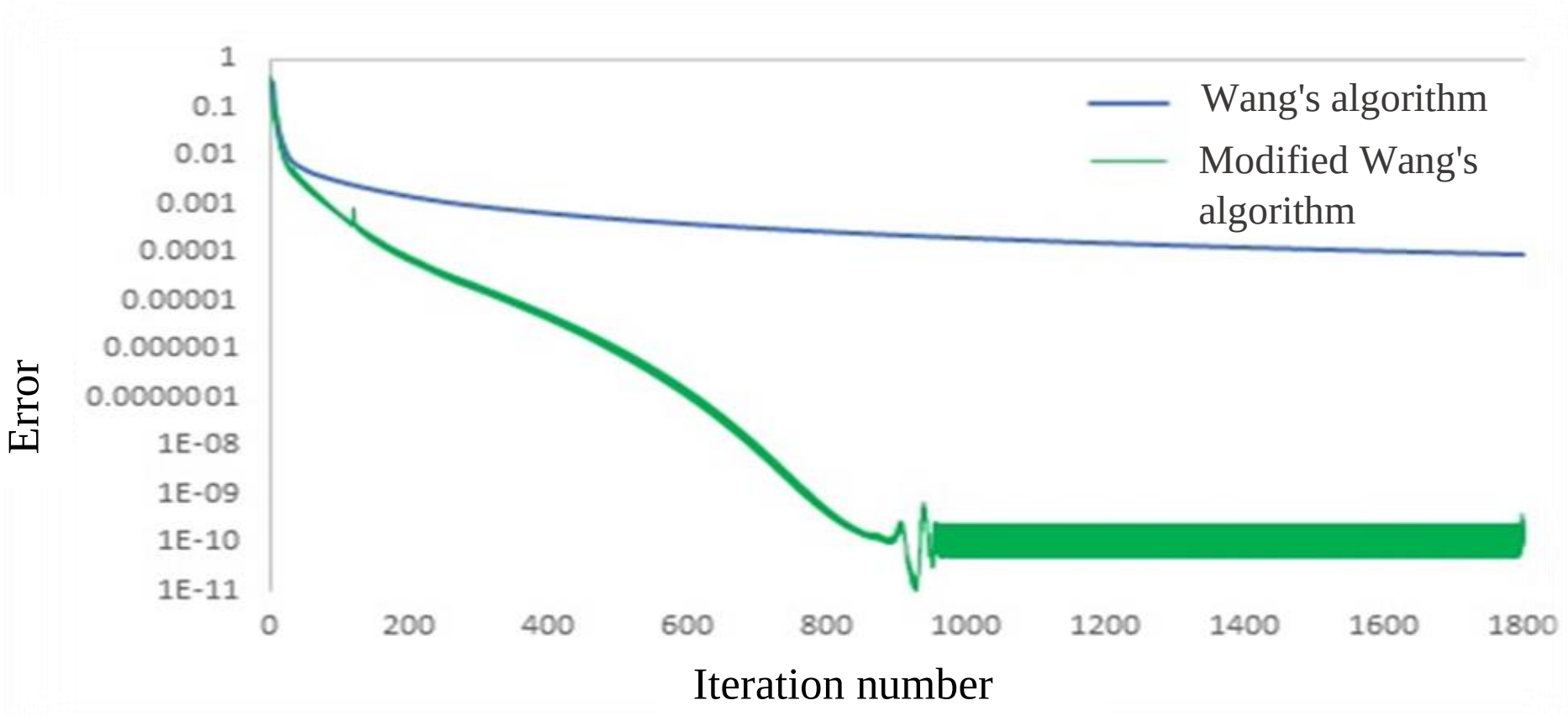

**Figure 2 Comparison of the convergence process of Wang and Modified Wang's algorithms in Nguyen network.**

**Sioux Falls network**

**Figure 3** shows the Sioux Falls network. Counting all paths in this network is impossible due to its large size. Therefore, in this section, only a limited number of paths will be considered for each algorithm.

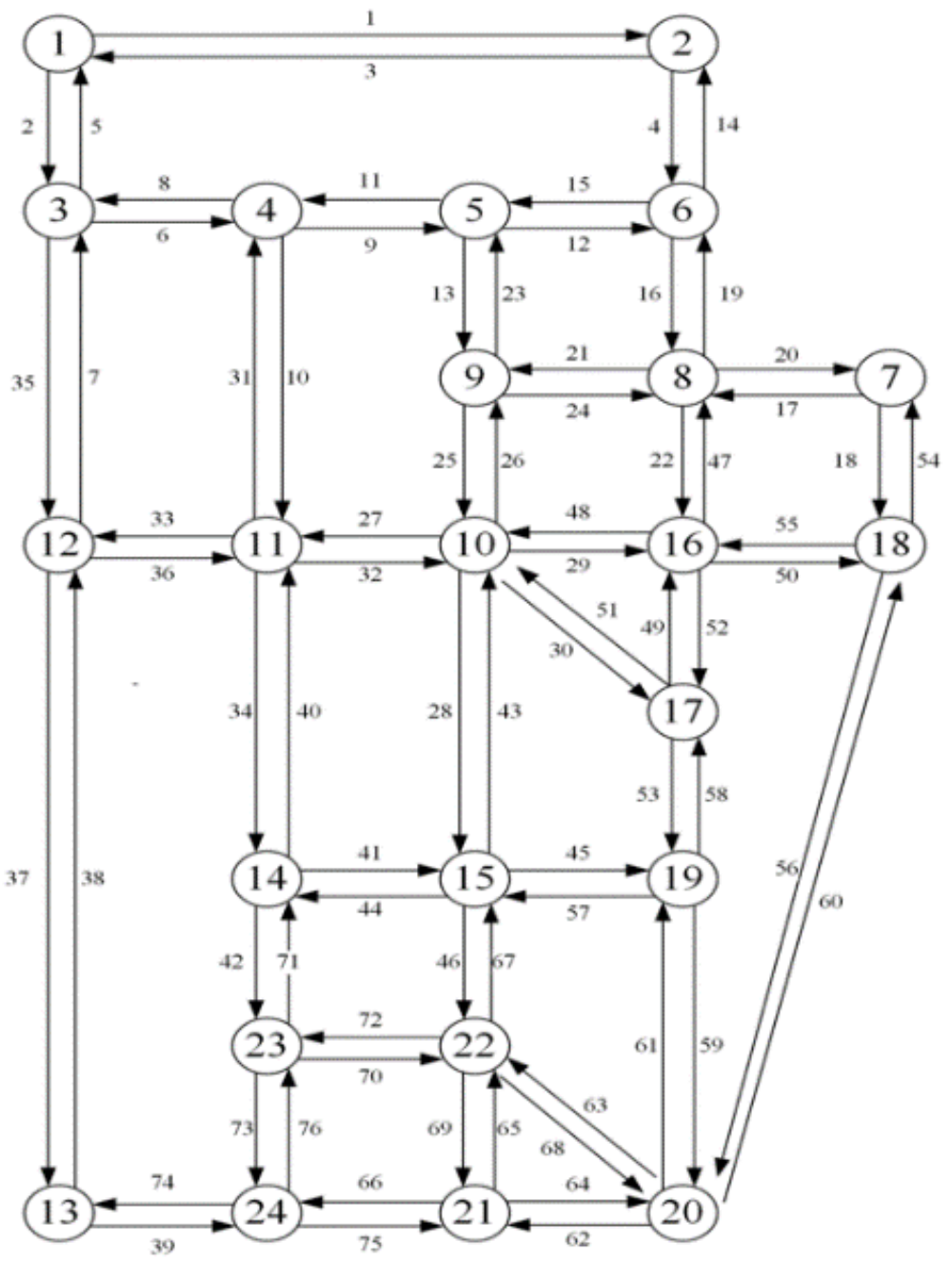

**Figure 3 Sioux falls network.**

**Implementation of two methods**

The Sioux Falls network has 528 OD pairs. **Table 2** shows the results of Wang's and Modified Wang's algorithms implemented in this network. For each OD, ten paths are considered. According to this table, the Modified Wang's algorithm achieves an error of 0.005, about 23 times faster than the Wang's algorithm.

**Table2 The results of two traffic assignment methods with ten routes for the Sioux Falls network**

| **solution method** | **Error** | **number of iterations** | **Execution time (seconds)** | **Total travel cost ($)** |
|---|---|---|---|---|
| Wang's algorithm | 0.005 | 2454 | 1949.72 | 1668808.96 |
| Modified Wang's algorithm | 0.005 | 125 | 87.95 | 1668888.96 |

**Figure 4** shows the convergence process of both Wang's and Modified Wang's algorithms for 1000 iterations so that the error changes in each iteration can be compared.

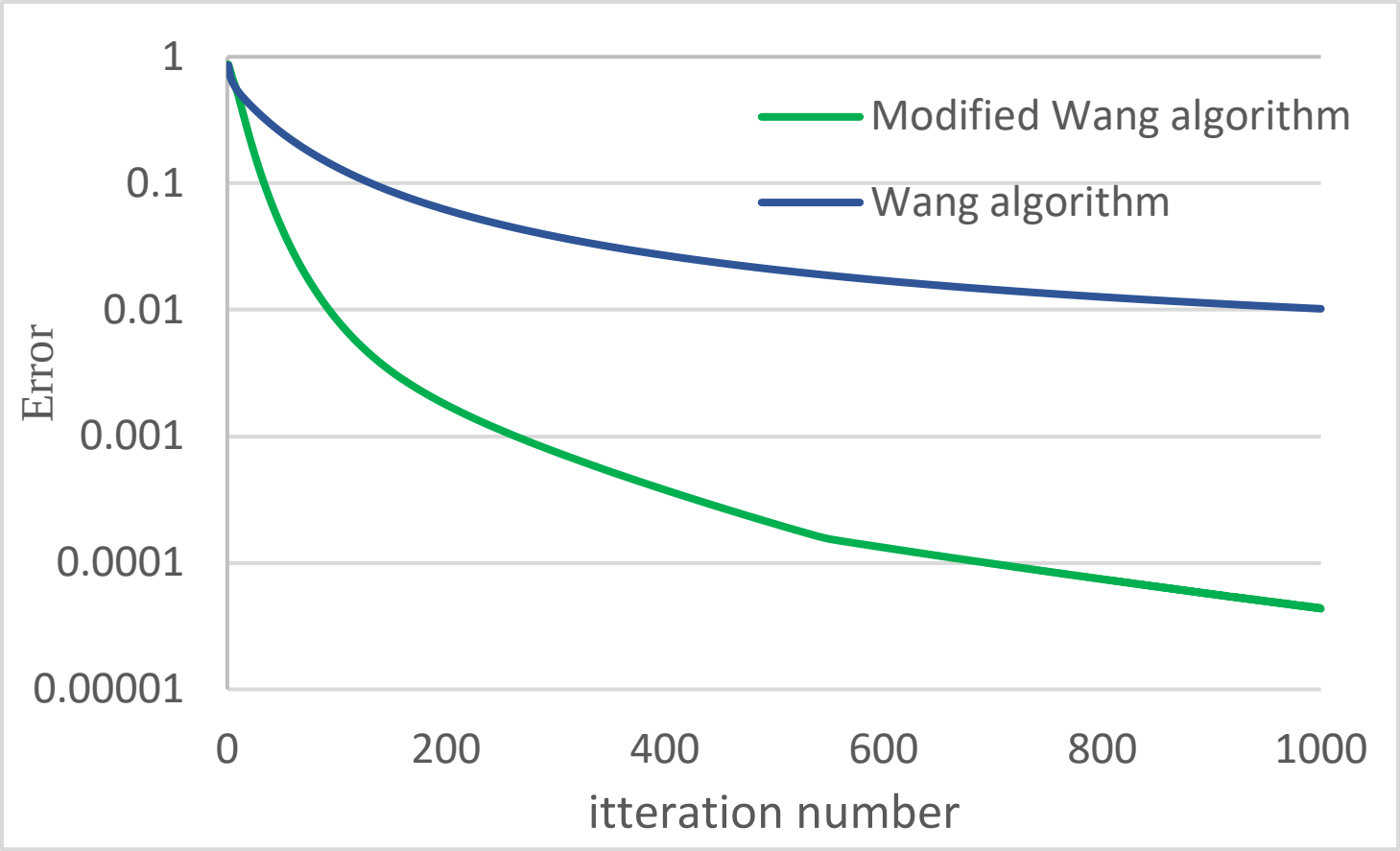

**Figure 4 A comparison between Wang's and Modified Wang's algorithms in the Sioux Falls network.**

**Implementation of path generation-assignment algorithm in Sioux Falls network**

In the traffic assignment process, all paths of the Sioux Falls network cannot be considered. Therefore, in all path-based traffic assignment algorithms, we must use a limited number of paths for each OD. This study will use the PGA algorithm due to its distinctive feature. This algorithm considers a limited number of paths and tries to choose the most attractive paths.

**Table 3** shows the results of the PGA algorithm with $k$ shortest paths for different values of $k$ for the Sioux Falls network. A comparison of the two flows can be made by computing the coefficient of determination and the flow deviation values. The amount of flow deviation ($dev$) is calculated according to the following equation:

$$dev_m^d(k) = \frac{\sum_a \left| x_{a,m}^{\mathrm{d}}(k) - x_{a,m}^{\mathrm{d}}(\bar{k}) \right|}{\sum_a x_{a,m}^{\mathrm{d}}(\bar{k})} \qquad \forall d, g, m \qquad (33)$$

According to **Equation 33**, $dev_m^d$ refers to the difference between the flow in the link resulting from the implementation of the $d$ solution method, considering $k$ path, and that from the implementation of the $d$ solution method, considering $\bar{k}$ path for the type $m$ vehicle. $x_{a,m}^{\mathrm{d}}$ respectively represent the flow in link $a$ for car type $m$ resulting from the implementation of $d$ solutions.

$dev_m^d$ has been calculated assuming 20 paths per OD for Sioux Falls; therefore, **Equation 33** considers $\bar{k}$ to be 20 for Sioux Falls. Based on **Table 3**, the PGA algorithm for $k = 20$ takes 35 times longer than for $k = 5$. For both types of vehicles, the coefficient of determination of the flow in the links is about 0.99 for $k = 20$ and $k = 5$. Therefore, their results are very similar.

**Figure 5** shows how the Modified Wang's algorithm's execution time changes for different paths for each OD. This figure shows that execution time is proportional to path number with a power of two, based on experimental data.

**TABLE 3 The results of the implementation of the PGA algorithm with k path for different values of k in the Sioux Falls network.**

| | ***k*=5** | ***k*=8** | ***k*=10** | ***k*=15** | ***k*=20** |
|---|---|---|---|---|---|
| The number of iterations of the k shortest path algorithm | 4 | 3 | 3 | 3 | 3 |
| Total number of flow change iterations | 273 | 397 | 493 | 737 | 1039 |
| Execution time (seconds) | 72.93 | 229.96 | 380.7 | 1129.46 | 2498.27 |
| Total travel cost ($) | 1543684 | 1628933 | 1667718 | 1716039 | 1741963 |
| Regular Vehicle $R^2$ | 0.9947 | 0.9992 | 0.9997 | 1.0000 | 1 |
| Autonomous Vehicle $R^2$ | 0.9957 | 0.9993 | 0.9995 | 0.9999 | 1 |
| RV $dev$ | 0.0591 | 0.0292 | 0.0067 | 0.0067 | 0 |
| AV $dev$ | 0.0168 | 0.0082 | 0.0022 | 0.0022 | 0 |

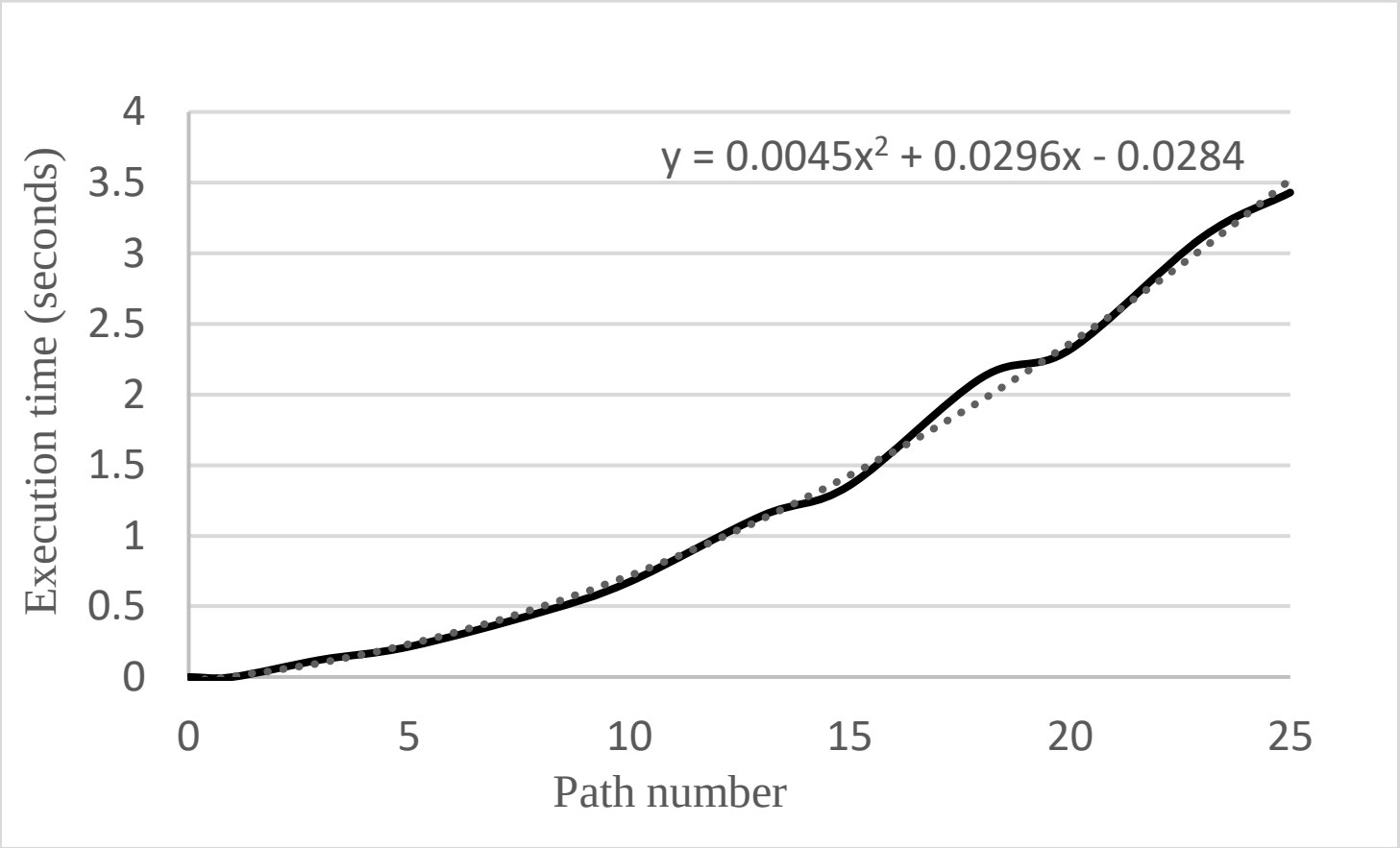

**Figure 5 Variations in the execution time of updating the flow of paths in a Modified Wang's algorithm iteration according to the number of paths in the Sioux Falls network.**

**Comparison of path generation-assignment algorithm with the Modified Wang's algorithm in Sioux Falls network**

In Modified Wang's algorithm, the $k$ shortest path is generated once based on the free travel cost of links and using Yen's algorithm. However, in the PGA algorithm, shortest path generation and traffic assignment are done alternately. As a result of this action, the flow deviation in the link is claimed to be reduced in comparison with the state in that all paths are considered. To better understand the results,

**Figures 6 and 7** show the Q-Q plots for the link resulting from both algorithms compared with the case where 20 paths are considered for each OD for both types of vehicles.

In **Figures 6 and 7**, RVs have a higher flow deviation than AVs for all values of $k$. This can be explained by the fact that in the UE state, only a few paths have flow, so by increasing the value of $k$, the paths with the flow will maintain their previous status.

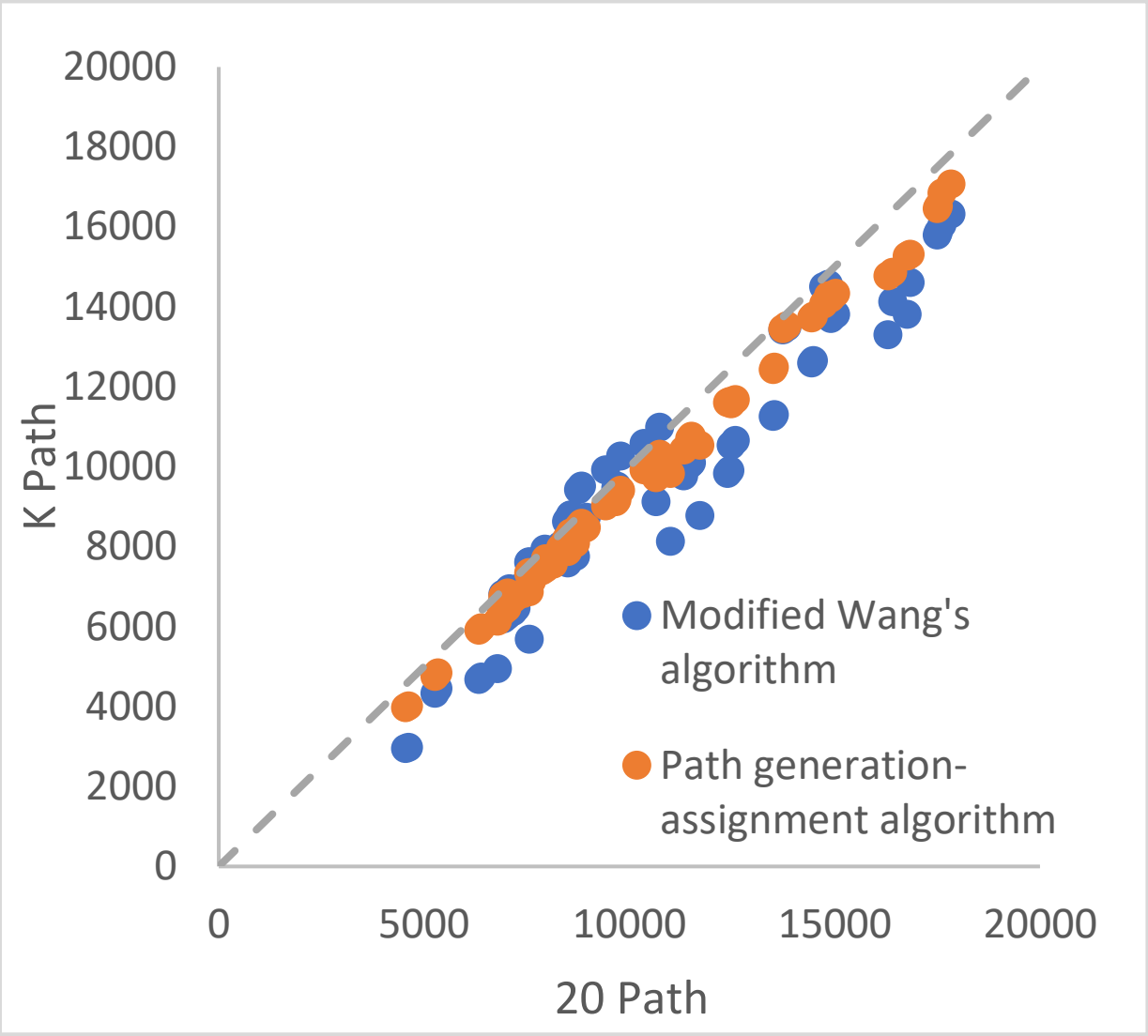

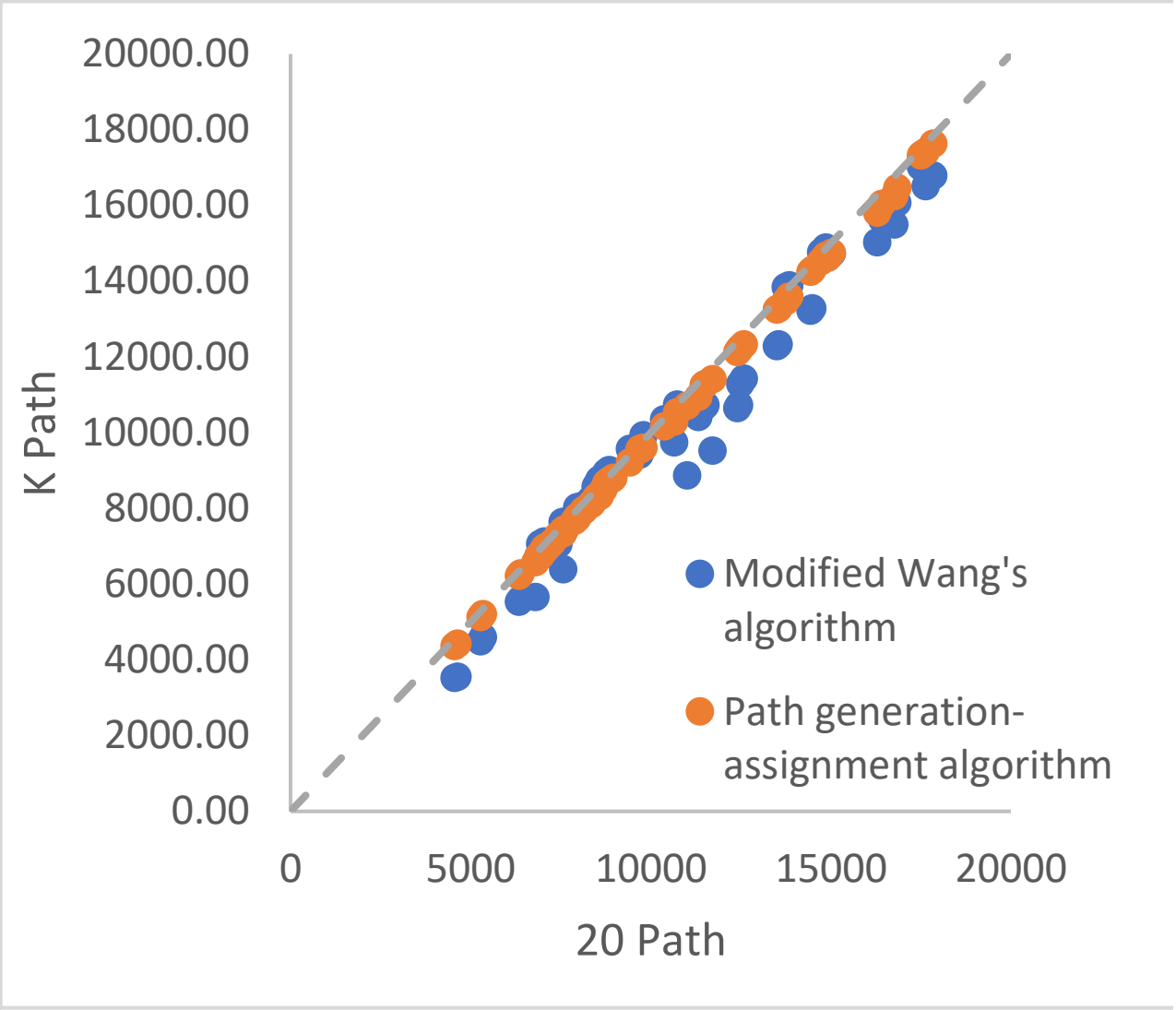

$(a)\ k = 5$ $(b)\ k = 10$

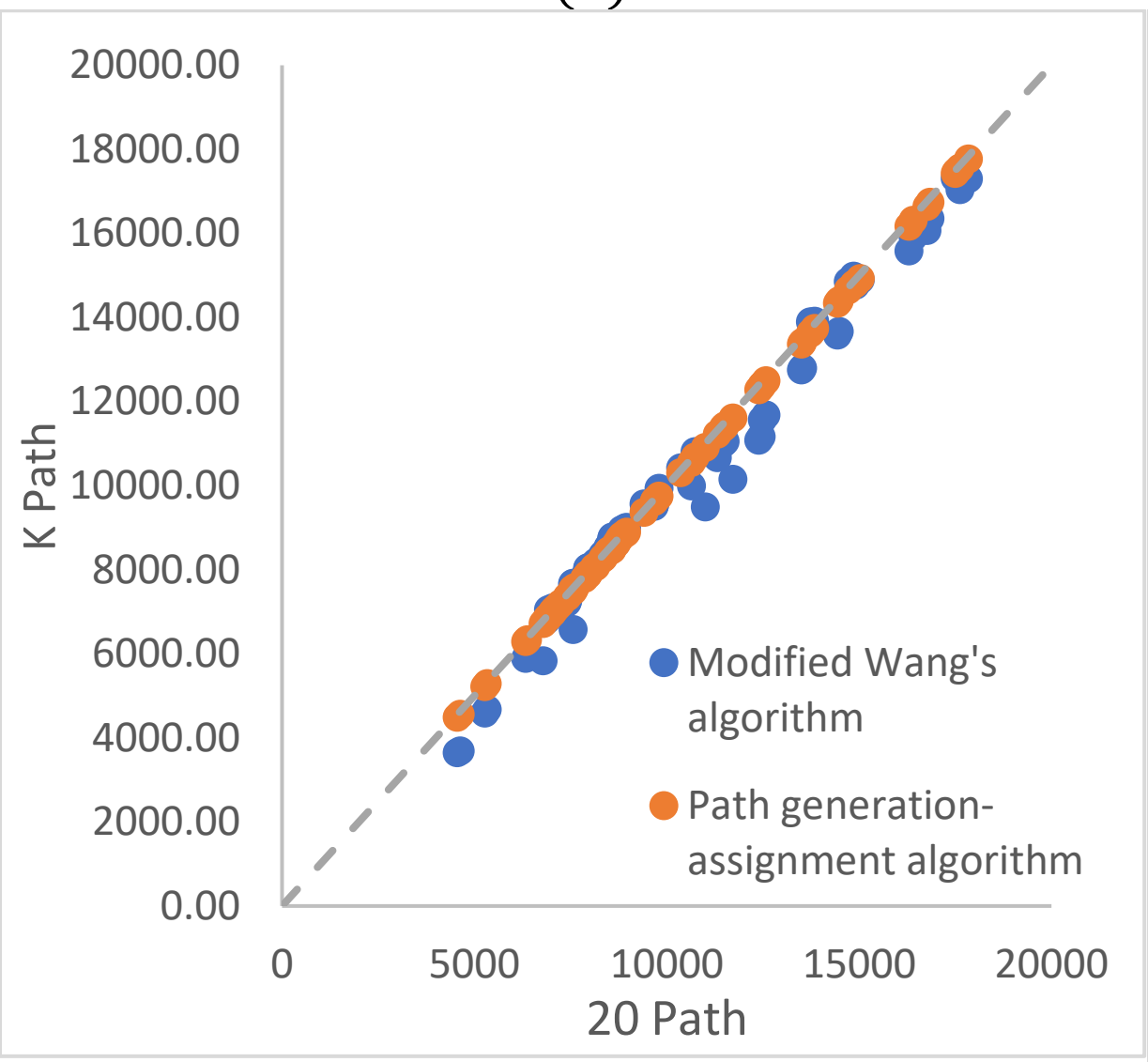

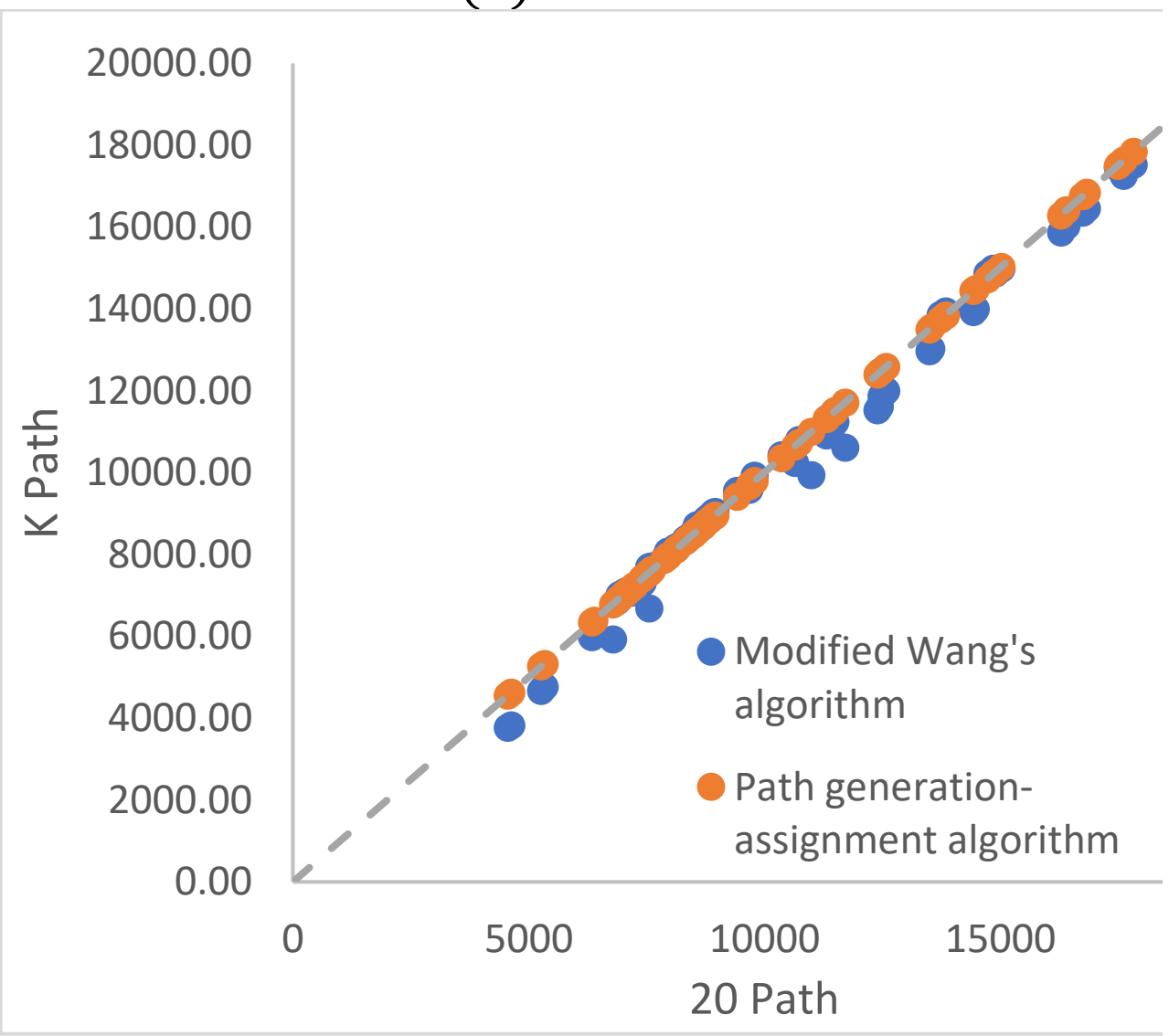

$(c)\ k = 15$ $(d)\ k = 20$

**Figure 6 Comparison of the flow in the link of RVs in the mode of using k paths with the mode of 20 paths for two Modified Wang's algorithms and PGA algorithms for different values of k in the Sioux Falls network.**

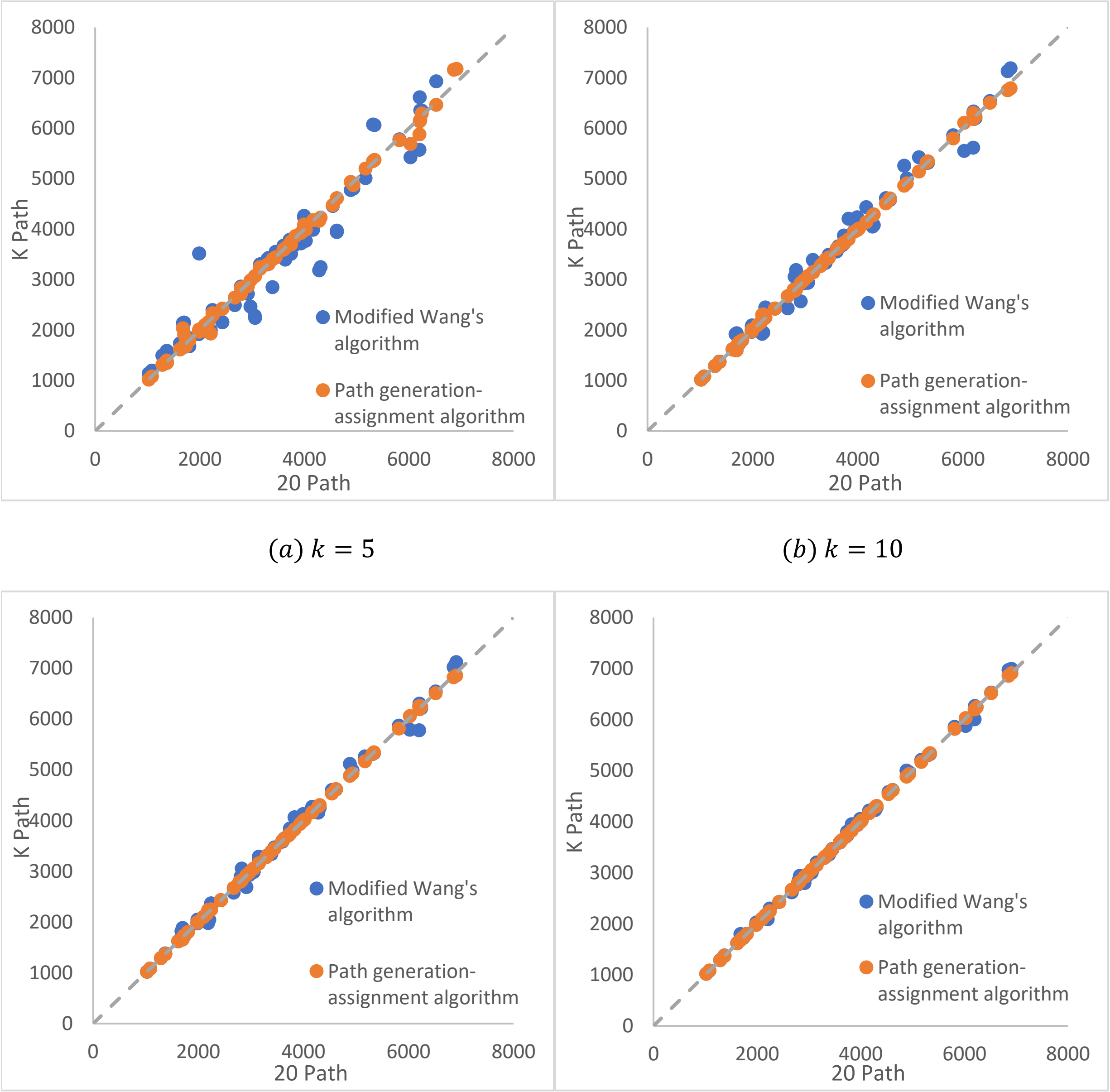

$(a)\ k = 5$ $(b)\ k = 10$

$(c)\ k = 15$ $(d)\ k = 20$

**Figure 7 Comparison of the flow of AVs in links in the mode of using k path with the mode of $k = 20$ for two Modified Wang's algorithms and PGA algorithms for different values of k in the Sioux Falls network**.

In **Figure 6**, the scatter of points around the line with a slope of 1 decreases as $k$ increases. Additionally, the diffusion of points related to the implementation of the PGA algorithm (orange points) is lower than that of points related to the Modified Wang's algorithm (blue points). Consequently, the flow in the RV link resulting from the implementation of the PGA is more similar to the case where 20 paths are considered for each OD.

**Figure 7** shows that the diffusion of flow points in the links of AVs decreases as k increases. According to this figure, the flow in the link of AVs resulting from the implementation of the PGA algorithm (orange points) has less scattering than the flow in the link of AVs resulting from the implementation of the Modified Wang's algorithm (blue points).

As shown in **Figures 6 and 7**, the PGA algorithm provides a better solution for the SMTA problem than the Modified Wang's algorithm in terms of flow deviation in the link.

## CONCLUSIONS

The article discussed the problem of SMTA. The transportation network was assumed to have two types of users. Each of these two classes has its own characteristics and preferences. The first class uses RVs and the second class uses AVs. Since RV users have incomplete knowledge of the network situation, the CNL model is used to model their behavior in the route choice stage. On the other hand, the users of AVs have full knowledge of the network situation due to the high technology of this category of vehicle, and therefore UE model is used for them in the route choice stage.

Three solution methods were presented for solving the problem. In Method 1, the SMTA problem is implemented as an NCP in GAMS. The second method is an iterative algorithm. This method decomposes the transportation network into OD units, and GAMS software solves an NCP for each OD iteration. In method 3, the Modified Wang's algorithm was introduced. The algorithm has increased convergence speed by altering some parameters used by Wang's algorithm in calculating descent direction and determining step size. Using Nguyen and Sioux Falls networks, we compared the convergence of this algorithm and Wang's algorithm.

Typically, path-based assignment algorithms require the path to be set before they can begin working. In some cases, such as Sioux Falls, where the dimensions of the network are relatively large, and it is not possible to determine all the network paths, the primary paths used by all three methods are generated using the k-shortest path algorithm. This research proposes a PGA algorithm to reduce the deviation of the final solution compared to a case where all paths are considered. Also, the comparison results of this algorithm and the Modified Wang's algorithm were shown for the Sioux Falls network, and the superiority of this algorithm compared to the Modified Wang's algorithm was emphasized.

## AUTHOR CONTRIBUTIONS

The authors confirm contribution to the paper as follows: study conception and design: SR. Mousavi, A. Yazdiani; data collection: SR. Mousavi, A. Yazdiani; analysis and interpretation of results: SR. Mousavi, A. Yazdiani; draft manuscript preparation: SR. Mousavi, A. Yazdiani. All authors reviewed the results and approved the final version of the manuscript.

**REFRENCES**